\documentclass[11pt]{article}
\usepackage{graphics}
\usepackage{graphicx}
\usepackage[all]{xy}
\usepackage{ifthen}
\usepackage{color}
\newfont{\bbb}{msbm10 scaled\magstephalf}
\newfont{\sbbb}{msbm7 scaled\magstephalf}
\def\C{\mbox{\bbb{C}}}
\def\Q{\mbox{\bbb{Q}}}
\def\R{\mbox{\bbb{R}}}
\def\Z{\mbox{\bbb{Z}}}
\def\cd{\C^d}
\def\rd{\R^d}
\def\rn{\R^n}
\def\zd{\Z^d}
\def\td{T^d}
\def\rndu{(\rn)^*}
\def\rddu{(\rd)^*}
\def\et1{e^{2\pi i\theta_1}}
\def\etd{e^{2\pi i\theta_d}}
\def\vz{\underline{z}}
\def\vu{\underline{u}}
\def\zjs{|z_j|^2}
\def\cl{{\cal C}_\lambda}
\def\D{\Delta}
\def\xd{X_1,\ldots,X_d}

\newtheorem{thm}{Theorem}[section]

\newtheorem{prop}[thm]{Proposition}
\newtheorem{remark}[thm]{Remark}

\newtheorem{example}[thm]{Example}
\newtheorem{ass}[thm]{Assumption}

\definecolor{red}{rgb}{.6,0,0}
\definecolor{green}{rgb}{0,.6,0}
\definecolor{darkgreen}{rgb}{0,0.3,0}
\definecolor{purple}{rgb}{0.5,0,0.5}
\definecolor{darkblue}{rgb}{0,0,0.7}
\definecolor{greenblue}{rgb}{0,0.4,0.5}
\newboolean{draft}
\newcommand{\cmt}[1]
{\ifthenelse {\boolean{draft}}
{{\sc \tiny \color{red} #1}}
{}}
\newcommand{\newbb}[1]
{\ifthenelse {\boolean{draft}}
{{\color{darkblue} #1}}
{#1}}
\newcommand{\newbbb}[1]
{\ifthenelse {\boolean{draft}}
{{\color{greenblue} #1}}
{#1}}
\newcommand{\inred}[1]
{\ifthenelse{\boolean{draft}}{{\color{red} #1}}{#1}}
\newcommand{\new}[1]
{\ifthenelse {\boolean{draft}}
{{\color{green} #1}}
{#1}}
\newcommand{\neww}[1]
{\ifthenelse {\boolean{draft}}
{{\color{darkgreen} #1}}
{#1}}
\newcommand{\newb}[1]
{\ifthenelse {\boolean{draft}}
{{\color{blue} #1}}
{#1}}
\newcommand{\del}[1]
{\ifthenelse {\boolean{draft}}
{{\color{magenta} #1}}
{}}
\setboolean{draft}{true}

\def\squareforqed{\hbox{\rlap{$\sqcap$}$\sqcup$}}
\def\qed{\ifmmode\else\unskip\quad\fi\squareforqed}

\newcommand{\proof}{\mbox{\bf Proof.\ \ }}

\newcounter{sect}

\title{\sc Nonrational Symplectic Toric Cuts}
\author{\sc Fiammetta Battaglia and Elisa Prato}
\date{}
\begin{document}
\maketitle
\begin{abstract} In this article we extend cutting and blowing up to the nonrational symplectic toric setting. 
This entails the possibility of cutting and blowing up for symplectic toric manifolds and orbifolds in nonrational directions.
\end{abstract}
{\small Mathematics Subject Classification 2010. Primary: 53D20
Secondary: 14M25}
\section*{Introduction}
\label{intro}
Symplectic blowing up was introduced by Gromov \cite{gromov} and later developed, from different points of view, by McDuff 
\cite{md0} and Guillemin--Sternberg \cite{gs1}. The symplectic cutting procedure, which was introduced by Lerman in \cite{lerman}, 
is an important generalization of the symplectic blowing up construction; it relies on the connection between blow--ups and symplectic 
reduction that was established by Guillemin--Sternberg in \cite{gs1}. 

We recall that, by the Atiyah, Guillemin--Sternberg convexity theorem \cite{a, gs}, the
image of the moment mapping for a Hamiltonian torus action on a
compact, connected symplectic manifold is a rational convex polytope, called moment polytope. In this context, cutting
the symplectic manifold corresponds to cutting the moment polytope with an appropriate rational hyperplane \cite[Remark~1.5]{lerman}.

In the special case of symplectic toric manifolds, that is, when the torus action is effective and the manifold has twice the dimension of the torus, 
the Delzant theorem \cite{d} establishes a one--to--one correspondence between the manifolds and the moment polytopes. In this case, the cut spaces
turn out to be the symplectic toric manifolds corresponding to the cut polytopes. Similar results, due to Lerman and Tolman \cite{lertol, lerman}, 
hold in the case of orbifolds. 

One of the remarkable features of the Delzant theorem is that it provides an explicit construction of the symplectic manifold from the polytope, 
following the same principle that allows to construct a complex toric variety from a fan. When extending the Delzant procedure to {\em any} simple polytope,
the lattice is replaced by a {\em quasilattice} and rationality is replaced by {\em quasirationality}.
The resulting spaces are typically not Hausdorff: they are modeled by quotients of smooth manifolds by the action of discrete
groups\footnote{By discrete we mean countable with the discrete topology.} \cite{pcras, p}.  We refer to them as {\em symplectic toric quasifolds}.

In this article we define symplectic cutting for symplectic toric quasifolds with respect to a quasirational cutting hyperplane and we prove that, 
as in the classical case, the cut spaces are the symplectic toric quasifolds that correspond to the cut polytopes. 
We remark that we allow cuts through vertices of the moment polytope: our only requirement is that the cut polytopes have maximal 
dimension and are simple. Interestingly, this can be done also for symplectic toric manifolds and orbifolds (see Remark~\ref{liscio}).

As a special case of symplectic cutting we obtain, as in the classical case, blowing up for symplectic toric quasifolds.

Finally, we show that our viewpoint allows us to make sense of cutting and blowing up even if the quasirationality 
condition on the hyperplane is dropped. In particular, this allows cutting and blowing up for symplectic toric manifolds and orbifolds in nonrational directions.

It will be convenient for us to prove our results in the more general setting of simple pointed polyhedra. A pointed polyhedron is a 
finite intersection of closed half--spaces that has at least one vertex; one can show that
it is a convex polytope if, and only if, it does not contain a ray (see \cite{ziegler}).

This article is structured as follows: in Section~\ref{delzant} we adapt to pointed polyhedra the generalized Delzant procedure given in \cite{p}; 
in Section~\ref{cuts} we define cutting for symplectic toric quasifolds and we prove our main result; in Section~\ref{cutkite},
as an example, we cut the Penrose kite in half; in Section~\ref{blowup} we discuss symplectic blowing up; finally,
in Section~\ref{arbitrary} we cut and blow up in arbitrary directions.
 
\section{The generalized Delzant procedure}
\label{delzant}
We refer the reader to \cite{p,kite} for the basic definitions and properties of quasifolds.
Following Ziegler \cite{ziegler}, we say that a subset $\D\subset\rndu$ is a {\em polyhedron} if it is a finite intersection of closed half--spaces. 
This implies, in particular, that $\D$ has at most a finite number of vertices. A polyhedron is said to be {\em pointed} if it has at least one
vertex. One easily sees that a pointed polyhedron is a convex polytope if, and only if, it does not contain a ray. 
As for a convex polytope, an $n$--dimensional pointed polyhedron is said to be {\em simple} if each of its vertices is contained in exactly $n$ facets.

A {\em quasilattice} $Q$ in $\rn$ is the $\Z$--span of a set of $\R$--spanning vectors $Y_1,\ldots,Y_d\in\rn$. Clearly, a quasilattice is an ordinary
lattice if, and only if,  it is generated by a basis of $\rn$. If $Q$ is a quasilattice in $\rn$, we say that the quasifold and abelian group 
$D^n=\rn/Q$ is a {\em quasitorus}.  As in the standard case, we will call the projection $\exp_{D^n}\,\colon\,\rn\longrightarrow D^n$ 
{\em exponential mapping}. Notice that, if $Q$ is an ordinary lattice, then $D^n$ is an ordinary torus.

Let now $\D\subset\rndu$ be an $n$--dimensional polyhedron having $d$ facets.
Then there exist $X_1,\ldots,X_d \in\rn$ and $\lambda_1,\ldots,\lambda_d\in\R$ such that
\begin{equation}\label{decomp}\D=\bigcap_{j=1}^d\{\;\mu\in\rndu\;|\;\langle\mu,X_j\rangle\geq\lambda_j\;\}.\end{equation}
Each of the vectors $X_1,\ldots,X_d$ is orthogonal to a facet of $\D$ and is inward--pointing.
We will say that $\D$ is {\em quasirational} with respect to the quasilattice $Q$ if the vectors $X_1,\ldots,X_d$ can be chosen in $Q$. 
If $\D$ is quasirational with respect to an ordinary lattice, then it is rational in the ordinary sense. Notice that any polyhedron is quasirational 
with respect to the quasilattice that is generated by a set of $d$ vectors, each of which is orthogonal to one of the $d$ different facets of $\D$.

We recall the generalized Delzant procedure for the nonrational setting \cite[Theorem~3.3]{p}. 
The variant presented here is a straightforward extension to the case of pointed polyhedra.
This allows a simultaneous treatment of polytopes and of pointed polyhedral cones, the latter being moment map 
images of equivariant local charts.  

\begin{thm}\label{thmp1}
Let $Q$ be a quasilattice in $\rn$ and let $\D\subset\rndu$ be an $n$--dimensional simple pointed polyhedron
that is quasirational with respect to $Q$. Assume that $d$ is the number of facets of $\D$ and
consider vectors $X_1,\ldots,X_d$ in $Q$ that satisfy (\ref{decomp}). For each $(\D, Q,\{\xd\})$,
there exists a $2n$--dimensional symplectic quasifold $(M,\omega)$
endowed with the effective Hamiltonian action of the quasitorus $D^n=\rn/Q$ such that, 
if $\Phi\,\colon M\rightarrow\rndu$ is the corresponding moment mapping, 
then $\Phi(M)=\Delta$. If $\Delta$ is a convex polytope, then $M$ is compact.
\end{thm}

\noindent \proof
Consider the space $\cd$
endowed with the standard symplectic form $\omega_0=\frac{1}{2\pi
i}\sum_{j=1}^d dz_j\wedge d\bar{z}_j$ and the standard action of the
torus $\td=\rd/\zd$:
$$
\begin{array}{cccccl}
\tau\,\colon& \td&\times&\cd&\longrightarrow& \cd\\
&((\et1,\ldots,\etd)&,&\vz)&\longmapsto&(\et1 z_1,\ldots, \etd z_d).
\end{array}
$$
This action is effective, Hamiltonian, and its moment mapping is given by
$$
\begin{array}{cccl}
J\,\colon&\cd&\longrightarrow &\rddu\\
&\vz&\longmapsto & \sum_{j=1}^d \zjs
e_j^*+\lambda,\quad\lambda\in\rddu \;\mbox{constant}.
\end{array}
$$
The mapping $J$ is proper and its image is the cone
$\cl=\lambda+{\cal C}_0$, where ${\cal C}_0$ denotes the positive orthant in the
space $\rddu$. Since $\Delta$ is pointed, the linear mapping
\begin{eqnarray*}
\pi\,\colon &\rd \longrightarrow \rn,\\
&e_j \longmapsto X_j
\end{eqnarray*}
is surjective.
Consider the $n$--dimensional quasitorus $D^n=\rn/Q$. Then the linear
mapping $\pi$ induces a quasitorus epimorphism $\Pi\,\colon\,\td
\longrightarrow D^n$. Define now $N$ to be the kernel of the mapping
$\Pi$ and choose $\lambda=\sum_{j=1}^d {\lambda_j} e_j^*$. Denote by
$i$ the Lie algebra inclusion $\mbox{Lie}(N)=\ker(\pi)\rightarrow\rd$ and
notice that $\Psi=i^*\circ J$ is a moment mapping for the induced
action of $N$ on $\cd$. The orbit space $M=\Psi^{-1}(0)/N$ is a 
symplectic quasifold by \cite[Theorem~3.1]{p}; its symplectic form $\omega$ is induced
from $\omega_0$. For the construction of an explicit atlas of $M$,
see the proof of \cite[Theorem~3.2]{cx}. 
The quasitorus $\td/N$ acts in a
Hamiltonian fashion on the symplectic quasifold
$M=\Psi^{-1}(0)/N$. If we identify the quasitori $D^n$ and $\td/N$
using the epimorphism $\Pi$, we get a Hamiltonian action of the
quasitorus $D^n$; the induced moment mapping is given by 
\begin{equation}\label{formula generale}\Phi=(\pi^*)^{-1}\circ J\end{equation}
and  
has image equal to
${(\pi^*)}^{-1}(\cl\cap\ker{i^*})=
{(\pi^*)}^{-1}(\cl\cap\mbox{im}\,\pi^*)= {(\pi^*)}^{-1}(\pi^*(\D))$,
which is exactly $\D$. This action is effective since the level set
$\Psi^{-1}(0)$ contains points of the form $\vz\in\cd$, $z_j\neq0$,
$j=1,\ldots,d$, where the $T^d$--action is free. Notice that
$\dim{M}=2d-2\dim{N}= 2d-2(d-n)=2n=2\dim{D^n}$. Finally we remark that 
$\Psi^{-1}(0)=J^{-1}(\pi^*(\Delta))$, therefore,
if $\Delta$ is a polytope, since $\pi^*$ is injective and $J$ is proper, $\Psi^{-1}(0)$ and thus $M$ are compact. 
\qed

This construction depends on two arbitrary
choices: the choice of the quasilattice $Q$ with respect to which
the pointed polyhedron is quasirational, and the choice of the
inward--pointing vectors $\xd$ in $Q$.
We will say that the quasifold $M$ with the effective Hamiltonian action of $D^n$ is the {\em symplectic toric quasifold} 
associated to $(\D, Q, \{\xd\})$. We recall that the case where $Q$ is a lattice and $\Delta$ is a rational simple convex polytope was treated 
by Lerman and Tolman \cite{lertol}, who allowed orbifold singularities and introduced the notion of 
{\em symplectic toric orbifold}. {\em Complex} toric quasifolds, on the other hand, were introduced and discussed in \cite{cx}.

Let now ${\widetilde Q}\subset Q$ be two quasilattices and assume that $\Delta$ is a simple pointed polyhedron that is quasirational with respect to 
$\widetilde Q$. Notice that $\Delta$ is quasirational with respect to $Q$ as well. 
Consider the quasitori ${\widetilde D}^n=\R^n/{\widetilde Q}$ and $D^n=\R^n/Q$;
if $\Gamma$ is the discrete quotient group $Q/{\widetilde Q}=\exp_{{\widetilde D}^n}(Q)$, then $D^n\simeq{\widetilde D}^n/\Gamma$.
\begin{prop}\label{sottoquasireticolo}If $\widetilde M$ and $M$ are the two quasifolds obtained by 
applying the generalized Delzant procedure to $(\D,{\widetilde Q}, \{\xd\})$ 
and $(\D, Q, \{\xd\})$, then  $M\simeq{\widetilde M}/\Gamma$. Moreover, if we denote by
$\widetilde \Phi$ and $\Phi$ the moment mappings corresponding to the actions of ${\widetilde D}^n$ and $D^n$ respectively, then
the following diagram commutes:
$$
\xymatrix{
{\widetilde M}\ar[dr]^{\widetilde \Phi}\ar[dd]&\\
&\qquad\Delta\subset(\R^n)^*\\
M\ar[ru]_{\Phi}&
}
$$
\end{prop}
\noindent\proof
Following the proof of Theorem~\ref{thmp1}, write $\widetilde{M}=\widetilde{\Psi}^{-1}(0)/{\widetilde N}$ and $M=\Psi^{-1}(0)/N$. 
Notice that $\widetilde{\Psi}=\Psi$ and that $N/\widetilde{N}\simeq\Gamma$. This implies that 
$M=\Psi^{-1}(0)/N\simeq(\widetilde{\Psi}^{-1}(0)/{\widetilde N})/(N/\widetilde{N})\simeq{\widetilde M}/\Gamma$. 
\qed
\section{Nonrational symplectic cutting}
\label{cuts}
In this section, we extend the cutting procedure \cite{lerman} to symplectic toric quasifolds. To do so, let us consider an $n$--dimensional simple 
pointed polyhedron $\Delta\subset\rndu$ that is quasirational with respect to a quasilattice $Q$.
Write $\D$ as in (\ref{decomp}), with $\xd\in Q$ and let $M$ be the symplectic toric quasifold obtained by
applying the generalized Delzant procedure to $(\D, Q, \{\xd\})$. We denote again by $\Phi$ the moment mapping for the corresponding $D^n$--action.
Let us now take $Y\in Q$, $\epsilon\in \R$ and let us 
cut the polyhedron $\D$ with the hyperplane 
$H(Y,\varepsilon)=\{\lambda\in\rndu\;|\;\langle \lambda,Y\rangle=\varepsilon\}$. 
\begin{ass}
The vector $Y$ and the parameter $\epsilon$ are chosen so that
$$\Delta_{\varepsilon^+}=\Delta\cap\{\lambda\in\rndu\;|\;\langle \lambda,Y\rangle\geq\varepsilon\}\quad\hbox{and}\quad
\Delta_{\varepsilon^-}=\Delta\cap\{\lambda\in\rndu\;|\;\langle\lambda,Y\rangle\leq\varepsilon\}$$
are $n$--dimensional simple pointed polyhedra.
\end{ass}
We define {\em symplectic cutting} for toric quasifolds as follows.
As a first step we consider the quasicircle subgroup 
$D^1\subset D^n$ generated by $Y$: 
$$D^1=\{\,\exp_{D^n}(tY)\in D^n\,|\,t\in \R\,\}.$$
The induced action of $D^1$ on $M$ is Hamiltonian with moment mapping given by the component $\Phi_Y=\langle \Phi,Y\rangle$.
Notice that $D^1\simeq \R/Q^1$, where $Q^1=\{\,l\in\R\;|\;lY\in Q\,\}$ is a quasilattice in $\R$ containing $\Z$. 
Let $\Lambda=\exp_{S^1}(Q^1)=\{e^{2\pi i t}\;|\;tY\in Q\}\subset S^1=\R/\Z$; then $D^1\simeq S^1/\Lambda$.
The quotient $\C/\Lambda$ is a symplectic quasifold.
\begin{remark}{\rm
Notice that if $Q$ is an ordinary lattice, and if $Y\in Q$ is not primitive, then $\Lambda$
is a finite, non--trivial subgroup of $S^1$ and $\C/\Lambda$ is an orbifold.}
\end{remark}
Consider now the two natural actions of
$D^1$ given by:
$$\begin{array}{ccc}
D^1\times\C/\Lambda&\longrightarrow&\C/\Lambda\\
(\exp_{D^1}(t),[z])&\longmapsto&[e^{\mp 2\pi i t}z]
\end{array}
$$
The product $M\times\C/\Lambda$ is a symplectic quasifold of dimension $2(n+1)$, endowed with two Hamiltonian effective
actions of the quasitorus $D^n\times D^1$. 
Let us now consider the two corresponding induced actions of $D^1$ on $M\times\C/\Lambda$; their moment maps are given by $\nu_{\mp} ((m,[w]))=
\Phi_Y(m)\mp|w|^2$. In agreement with the smooth setting \cite{lerman}, 
we call {\em symplectic cutting} the operation that produces the two orbit spaces 
$$\overline{M}_{\Phi_Y\geq \varepsilon}=\nu_{-}^{-1}(\varepsilon)/D^1\quad\hbox{and}\quad\overline{M}_{\Phi_Y\leq \varepsilon}=
\nu_{+}^{-1}(\varepsilon)/D^1.$$
The action of $D^n$ on the first factor of $M\times\C/\Lambda$ commutes with the $D^1$--action and induces an effective $D^n$--action on both.
We will describe the symplectic toric quasifold structure of the cut spaces explicitly in Theorem~\ref{teorema} below. 
The first step is to apply the generalized Delzant procedure to $\Delta_{\varepsilon^+}$ and $\Delta_{\varepsilon^-}$. 
As normal vectors for $\Delta_{\varepsilon^+}$  we consider the appropriate subset of $X_1,\ldots,X_d$ plus the vector $Y$.
We argue similarly for $\Delta_{\varepsilon^-}$, with the vector $-Y$ instead. We thus obtain symplectic toric quasifolds
$M_{\varepsilon^+}$ and $M_{\varepsilon^-}$. Then we have
\begin{thm}\label{teorema} 
There exist $D^n$--equivariant homeomorphisms
$$f_+\,\colon\, M_{\varepsilon^+}\rightarrow\overline{M}_{\Phi_Y\geq \varepsilon}$$
$$f_-\,\colon\, M_{\varepsilon^-}\rightarrow\overline{M}_{\Phi_Y\leq \varepsilon}.$$
Moreover, ${M}_{\Phi_Y\geq \varepsilon}$ and ${M}_{\Phi_Y\leq \varepsilon}$ naturally inherit symplectic toric quasifold structures from
$M_{\varepsilon^+}$ and $M_{\varepsilon^-}$ respectively.
\end{thm}
\proof
We prove the theorem for $M_{\varepsilon^+}$; the other case can be treated in the same way.
Let $\upsilon$ be a vertex in $\Delta_{\varepsilon^+}$ that does not lie in the hyperplane
$H(Y,\varepsilon)$.
After a possible reordering of the facets, we have $\upsilon=\cap_{j=1}^{n}\{\lambda\in\rndu\;|\;\langle \lambda,X_j\rangle= \lambda_j\}$ and
$$\Delta_{\varepsilon^+}=\left(\cap_{j=1}^{l}\{\lambda\in\rndu\;|\;\langle \lambda,X_j\rangle\geq \lambda_j\}\right)\cap\{\lambda\in\rndu
\;|\;\langle \lambda,Y\rangle\geq\varepsilon\},$$
where $l+1$ is the number of facets of $\Delta_{\varepsilon^+}$ and $n\leq l \leq d$.

Before we define the mapping $f_+\colon M_{\varepsilon^+}\rightarrow\overline{M}_{\Phi_Y\geq \varepsilon}$, let us give an explicit description of 
$M$ and $M_{\varepsilon^+}$. Following the proof of Theorem~\ref{thmp1}, let us consider the projections
\begin{equation}\begin{array}{ccccc}
\pi&\colon&\R^d&\longrightarrow&\R^n\\
&&e_j&\longmapsto&X_j
\end{array}
\end{equation}
and
\begin{equation}\begin{array}{cccccc}
\pi_{\epsilon^+}&\colon&\R^{l+1}&\longrightarrow&\R^n&\\
&&e_j&\longmapsto&X_j&j=1,\ldots,l\\
&&e_{l+1}&\longmapsto&Y&
\end{array}
\end{equation}
Consider $N=\{\exp_{T^d}(X)\in T^d\;|\;\pi(X)\in Q\}$ , $N_{\epsilon^+}=\{\exp_{T^{l+1}}(X)\in T^{l+1}\;|\:\pi_{\epsilon^+}(X)\in Q\}$, and the corresponding
moment mappings $\Psi$ and $\Psi_{\epsilon^+}$. We have that $M=\Psi^{-1}(0)/N$ and
$ M_{\varepsilon^+}=\Psi_{\epsilon^+}^{-1}(0)/N_{\epsilon^+}$. Let us then compute $N, N_{\epsilon^+}, \Psi$ and $\Psi_{\epsilon^+}$.
The $n\times d$ matrix associated to the projection
$\pi$ with respect to the standard basis of $\R^d$  and the basis $\{X_1,\ldots,X_n\}$ of $\R^n$
can be written as $(I_n,A)$.
Similarly, the $n\times(l+1)$ matrix associated to the projection
$\pi_{\epsilon^+}$ with respect to the standard basis of $\R^{l+1}$  and the basis $\{X_1,\ldots,X_n\}$ of $\R^n$
can be written as $(I_n,A^{+})$.
It will be convenient for us to assume from now on that $n<l<d$;
the cases $l=n$ and $l=d$ are simpler and are left to the reader. 
Then, $A^{+}=(A^{l},\underline{b})$, where 
$A^{l}$ is the $n\times(l-n)$ matrix obtained from $A$ by deleting its last $d-l$
columns and $\underline{b}$ is the column $(b_1,\ldots,b_n)^t$, where $Y=\sum_{j=1}^n b_jX_j$. 
Let us now choose $\{X_{d+1},\ldots,X_q\}$ in $\R^n$ such that $\{X_{1},\ldots,X_{d},X_{d+1},\ldots,X_q\}$
is a set of generators of the quasilattice $Q$.
Let $C$ be the $n\times q$ matrix whose columns are the components of the 
generators $X_j$ with respect to the basis $\{X_1,\ldots,X_n\}$. 
Therefore $C=(I_n,A,P)$, for some $n\times(q-d)$ matrix $P$.
Notice now that any element of the quasilattice $Q$, $\sum_{i=1}^qm_iX_i$, $\underline{m}=(m_1,\ldots,m_q)\in \Z^q$, 
can be written in terms of the basis $\{X_1,\ldots,X_n\}$ as
$\sum_{j=1}^n (C_j\cdot\underline{m}) X_j$, where $C_j\cdot\underline{m}$ denotes the standard scalar product of the $j$-th row of $C$ with $\underline{m}$.
We are now ready to compute the group $N$. 
The identity $\pi(x_{1},\ldots,x_d)=\sum_{i=1}^qm_iX_i$ becomes
$\sum_{j=1}^{n}\left(x_j-C_j\cdot\underline{m}+A_j\cdot\underline{x}\right)X_j=0$,
where $\underline{x}=(x_{n+1},\ldots,x_d)$. Therefore N equals
\begin{equation}\label{enne}
\exp_{T^d}\{(C_1\cdot\underline{m}- A_1\cdot\underline{x},
\ldots,C_n\cdot\underline{m}-A_n\cdot\underline{x},
\underline{x})\in\R^d\;|\;
\underline{m}\in\Z^q\,,\underline{x}\in\R^{d-n}\}.
\end{equation}
Similarly one shows that the group $N_{\epsilon^+}$ equals
\begin{equation}\label{ennepiu}
\exp_{T^{l+1}}\{(C_1\cdot\underline{m}-A^{+}_1\cdot\underline{x},\ldots,C_n\cdot\underline{m}-A^{+}_n\cdot\underline{x},
\underline{x})\in\R^{l+1}\;|\;
\underline{m}\in\Z^q\,,\underline{x}\in\R^{l+1-n}\}.
\end{equation}
where $\underline{x}=(x_{n+1},\ldots,x_{l+1})$.
Let us now compute the moment mappings $\Psi$ and $\Psi_{\epsilon^+}$. The vectors
$$V_j=(a_{1j},\ldots,a_{nj},0,\ldots,0,-1,0,\ldots,0),\quad j=n+1,\ldots,d$$
(the $1$ here appearing at the $j$--th place) form a basis for the Lie algebra of $N$. Similarly, the vectors
$$
\left\{\begin{array}{ll}
V^+_j&=(a_{1j},\ldots,a_{nj},0,\ldots,0,-1,0,\ldots,0),\quad j=n+1,\ldots,l\\
V^+_{l+1}&=(b_1,\ldots,b_n,0,\ldots,0,-1)
\end{array}\right.
$$
form a basis for the Lie algebra of $N_{\epsilon^+}$.
Let $\{\alpha_{n+1},\ldots,\alpha_{d}\}$  and  
$\{\alpha^+_{n+1},\ldots,\alpha^+_{l+1}\}$ be the
respective dual bases.
A straightforward computation yields:
$$\Psi(\vz)=\sum_{k=n+1}^d\left(
\sum_{j=1}^na_{jk}(|z_j|^2+\lambda_j)-|z_k|^2-\lambda_k\right)\alpha_k
$$
and
$$
\begin{array}{lll}
\Psi_{\epsilon^+}(\vz)&=&\sum_{k=n+1}^{l}\left(
\sum_{j=1}^na_{jk}(|z_j|^2+\lambda_j)-|z_k|^2-\lambda_k\right)\alpha^+_k+\\&+&
\left(\sum_{j=1}^nb_j(|z_j|^2+\lambda_j)-|z_{l+1}|^2-\epsilon\right)
\alpha^+_{l+1}.
\end{array}
$$
From here one derives explicit expressions for $M$ and $M_{\varepsilon^+}$.
Consider now the mapping
\begin{equation}\label{mappaeffe}
\begin{array}{ccccc}
\widetilde{f_+}&\colon&\Psi_{\varepsilon^+}^{-1}(0)&\longrightarrow&\overline{M}_{\Phi_Y\geq \varepsilon}\\
&&(u_1,\ldots,u_{l+1})&\longmapsto&\left[[u_1:\cdots:u_l:w_{l+1}:\cdots:w_{d}]:[u_{l+1}]\right]
\end{array}
\end{equation}
where
$$w_k=\sqrt{\sum_{j=1}^na_{jk}(|u_j|^2+\lambda_j)-\lambda_k}.$$
In order to check that the induced mapping $f_+\colon M_{\varepsilon^+}\rightarrow\overline{M}_{\Phi_Y\geq \varepsilon}$
is a well defined $D^n$--equivariant homeomorphism, it is necessary to write explicitly the action of $D^1$ on $M$. 
The condition
$\Pi(\exp_{T^d}(x_1,\ldots,x_d))=\exp_{D^n}(tY)$ yields, for $\underline{x}=(x_{n+1},\ldots,x_d)$,
$$\sum_{j=1}^{n}(x_j+A_j\cdot\underline{x})X_j\simeq tY(\hbox{mod}Q);$$
therefore there exists $\underline{m}\in\Z^q$ such that
$$\sum_{j=1}^{n}\left(x_j-C_j\cdot\underline{m}+
A_j\cdot\underline{x}-tb_j\right)X_j=0.$$
Since $\exp_{T^d}(C_1\cdot\underline{m}-
A_1\cdot\underline{x},\ldots,
C_n\cdot\underline{m}-
A_n\cdot\underline{x}
,\underline{x})\in N$,
the action $D^1\times M\longrightarrow M$ is given by
$$
\exp_{D^n}(t Y)\cdot[z_1:\cdots:z_d]=\left[e^{2\pi itb_1}z_1:\cdots:e^{2\pi itb_n}z_n:z_{n+1}:\cdots:z_d\right].
$$
Consistently, by (\ref{formula generale}), the moment mapping $\nu_{-}$ with respect to the $D^1$--action on $M\times\C/\Lambda$ has the
following form:
$$\begin{array}{ccccccc}
\nu_{-}&\colon&M&\times& \C/\Lambda &\longrightarrow&\R^*\\
&&([z_1:\cdots:z_d]&,&[w])&\longmapsto& b_1(|z_1|^2+\lambda_1)+\cdots+b_n(|z_n|^2+\lambda_n)- |w|^2 .
\end{array}
$$
We now show that $f_+$ is injective; showing that it is well--defined and surjective 
rests on similar arguments.
Consider two points
$[u_1:\cdots:u_{l+1}]$, $[u'_1:\cdots: u'_{l+1}]$ in $M_{\epsilon^+}$
such that $f_+([u_1:\cdots:u_{l+1}])=f_+([u'_1:\cdots:u'_{l+1}])$.
Then
$$
\begin{array}{l}
[[u_1:\cdots:u_{l}:w_{l+1}:\cdots:w_d]:[u_{l+1}]]=[[u'_1:\cdots:u'_{l}:w'_{l+1}:\cdots:w'_d]:[u'_{l+1}]].
\end{array}
$$
Therefore there exists $t\in\R$ such that
$$
\begin{array}{l}
([e^{2\pi itb_1}u_1:\cdots:e^{2\pi itb_n}u_n:u_{n+1}:\cdots:u_l:w_{l+1}:\cdots:w_d],[e^{-2\pi it}u_{l+1}])=\\
([u'_1:\cdots:u'_{l}:w'_{l+1}:\cdots:w'_d],[u'_{l+1}]).
\end{array}
$$
This equation implies that there exists $\underline{m}\in\Z^q,\underline{x}\in\R^{d-n},\theta\in Q^1$ 
such that:
$$
\left\{\begin{array}{ll}
e^{2\pi i(C_j\cdot\underline{m}-A_j\cdot\underline{x}+tb_j)}u_j=u'_j&\;\;j=1,\ldots,n\\
e^{2\pi i x_{j}}u_{j}=u'_{j}&\;\;j=n+1,\ldots,l\\
e^{2\pi i x_{j}}w_{j}=w'_{j}&\;\;j=l+1,\ldots d\\
e^{2\pi i(\theta-t)}u_{l+1}=u'_{l+1}.&
\end{array}\right.
$$
This implies that $(x_{l+1},\ldots,x_d)\in\Z^{d-l-1}$, therefore 
$$\underline{n}=(m_1,\ldots,m_l,m_{l+1}+x_{l+1},\ldots,m_d+x_d,m_{d+1},\ldots,m_q)\in\Z^q$$ 
and $C_j\cdot\underline{m}-A_j\cdot\underline{x}+tb_j=C_j\cdot\underline{n}-A^l_j\cdot(x_{n+1},\ldots,x_l)+tb_j$.
Now consider $s=\theta-t$ and notice that, since $\theta Y\in Q$, there exists $\underline{k}\in\Z^q$ such that
$\theta b_j=C_j\cdot\underline{k},\;j=1,\ldots,n$. Thus we obtain
$$C_j\cdot\underline{m}-A_j\cdot\underline{x}+tb_j=C_j\cdot(\underline{n}+\underline{k})-A^+_j\cdot(x_{n+1},\ldots,x_l,s).
$$
Therefore
$$
\left\{\begin{array}{ll}
e^{2\pi i(C_j\cdot(\underline{n}+\underline{k})-A^+_j\cdot(x_{n+1},\ldots,x_l,s))}u_j=u'_j&\;\;j=1,\ldots,n\\
e^{2\pi i x_{j}}u_{j}=u'_{j}&\;\;j=n+1,\ldots,l\\
w_{j}=w'_{j}&\;\;j=l+1,\ldots d\\
e^{2\pi i s}u_{l+1}=u'_{l+1}.&
\end{array}\right.
$$
By (\ref{ennepiu}), we conclude that
$[u_1:\cdots:u_{l+1}]=[u'_1:\cdots:u'_{l+1}]$ and thus $f_+$ is injective.
It is straightforward to prove that the map $f_+$ is continuous. In order to prove that $f_+$ is open, we consider the 
quasifold atlas for the toric quasifold $M_{\varepsilon^+}$ given in the proof of \cite[Theorem~3.2]{cx}; its 
charts, $(U_{\mu},\psi_\mu,{\tilde U}/\Gamma_\mu)$, are indexed by the vertices $\mu$ of $\Delta_{\varepsilon^+}$. 
Consider one such vertex $\mu$ and the corresponding chart 
$${\widetilde U}_{\mu}/\Gamma_{\mu}\stackrel{\psi_{\mu}}{\longrightarrow}U_{\mu}\subset M_{\varepsilon^+}.$$
It is sufficient to prove that for each $\mu$ the mapping $f_+\circ\psi_{\mu}$ is open.
Consider $[v_1:\cdots:v_n]\in{\widetilde U}_{\mu}/\Gamma_{\mu}$; we have 
$\psi_{\mu}([v_1:\cdots:v_n])=[u_1(\underline{v}):\cdots:u_{l+1}(\underline{v})],$ where
the  functions $u_j(\underline{v})$ are given by
$$u_j(\underline{v})=\left\{\begin{array}{ll}
v_j &j=1,\ldots, n \\
\sqrt{\sum_{h=1}^{n}a_{jh}(|v_h|^2+\lambda_h)-\lambda_j}\quad &j=n+1,\ldots,l \\
\sqrt{\sum_{h=1}^{n}b_h(|v_h|^2+\lambda_h)-\varepsilon}\quad &j=l+1.
\end{array}\right.
$$
Therefore
$(f_+\circ\psi_{\mu})([\underline{v}])=[[u_1(\underline{v}):\cdots:u_l(\underline{v}),w_{l+1}(\underline{v}):\cdots:w_d(\underline{v})]:[u_{l+1}(\underline{v})]]$.
Let $W$ be an open subset in ${\widetilde U}_{\mu}/\Gamma_{\mu}$ and let $\widetilde W$ the $\Gamma_{\mu}$--saturated
open subset of ${\widetilde U}_\mu$ projecting to $W$. Consider the natural projections 
$p_1\,\colon \Psi^{-1}(0)\times\C\rightarrow(\Psi^{-1}(0)/N)\times(\C/\Lambda)$ and 
$p_2\,\colon\nu_-^{-1}(\varepsilon)\rightarrow\nu_-^{-1}(\varepsilon)/D^1$.
Then the set $f_+\circ\psi_{\mu}(\widetilde{W})$ is open in $\nu_-^{-1}(\varepsilon)/D^1$, if and only if the set
${\widehat W}=(p_2\circ p_1)^{-1}\left(f_+\circ \psi_{\mu}(\widetilde{W})\right)$ is open in
$(\nu_-\circ p_1)^{-1}(\varepsilon)$.
On the other hand, 
$\widehat W$ is the set of points
$(z_1,\ldots,z_d,w)\in\Psi^{-1}(0)\times\C$ such that
$$
\left\{\begin{array}{ll}
z_j=e^{2\pi i(C_j\cdot\underline{m}-A_j\cdot\underline{x}+tb_j)}u_j(\underline{v})&\;\;j=1,\ldots,n\\
z_j=e^{2\pi i x_{j}}u_{j}(\underline{v})&\;\;j=n+1,\ldots,l\\
z_j=e^{2\pi i x_{j}}w_{j}(\underline{v})&\;\;j=l+1,\ldots, d\\
w=e^{2\pi i(\theta-t)}u_{l+1}(\underline{v}).&
\end{array}\right.
$$
with $\underline{v}\in{\tilde W}$, $(x_{n+1},\ldots,x_d)\in\R^{d-n}$, $\underline{m}\in\Z^q$, $t\in\R$ and $\theta\in Q^1$.
Now observe that, arguing as in the proof of the injectivity of $f_+$, we have
$$C_j\cdot\underline{m}-A_j\cdot\underline{x}+tb_j=C_j\cdot(\underline{m}+\underline{k})-A_j^+\cdot(x_{n+1},\ldots,x_l,s)
-A_j(0,\ldots,0,x_{l+1},\ldots,x_d).$$
Setting $\underline{h}=\underline{m}+\underline{k}\in\Z^q$, notice that, by (\ref{ennepiu}),
$$
\left\{\begin{array}{ll}
e^{2\pi i(C_j\cdot\underline{h}-A_j^+\cdot(x_{n+1},\ldots,x_l,s))}u_j&\;\;j=1,\ldots,n\\
e^{2\pi i x_{j}}u_{j}&\;\;j=n+1,\ldots,l\\
e^{2\pi i s}u_{l+1}.&
\end{array}\right.
$$
describes the action of $N_{\varepsilon^+}$ on 
$[u_1:\cdots:u_{l+1}]\in M_{\varepsilon^+}$.
Now recall from (ii) in \cite[Lemma~2.3]{cx} that $N_{\varepsilon^+}=\Gamma_{\mu}\exp(\hbox{Lie}(N_{\varepsilon^+}))$, 
where the action of $\exp(\hbox{Lie}(N_{\varepsilon^+}))$ is obtained by setting $\underline{h}=0$ in the expression above.
Since $\widetilde W$ is $\Gamma_{\mu}$-invariant, the set $\widehat W$ is equal to the set of points
$(z_1,\ldots,z_d,w)\in\Psi^{-1}(0)\times\C$ such that
\begin{equation}\label{ultima}
\left\{\begin{array}{ll}
z_j=e^{2\pi i(-A_j^+\cdot(x_{n+1},\ldots,x_l,s)
-A_j(0,\ldots,0,x_{l+1},\ldots,x_d))}u_j(\underline{v})&\;\;j=1,\ldots,n\\
z_j=e^{2\pi i x_{j}}u_{j}(\underline{v})&\;\;j=n+1,\ldots,l\\
z_j=e^{2\pi i x_{j}}w_{j}(\underline{v})&\;\;j=l+1,\ldots, d\\
w=e^{2\pi i s}u_{l+1}(\underline{v}).&
\end{array}\right.
\end{equation}
where $\underline{v}\in{\widetilde W}$, $(x_{n+1},\ldots,x_d)\in\R^{d-n}$ and $s\in\R$.
Therefore $\widehat W$ is open since it is the image of ${\widetilde W}\times\R^{d-n}\times\R$
via the open mapping $ {\tilde U}_{\mu}\times\R^{d-n}\times\R\rightarrow (\nu_-\circ p_1)^{-1}(\varepsilon)$ that sends
$(\underline{v},x_{n+1},\ldots,x_d,s)$ to $(z_1,\ldots,z_d,w)$ given by (\ref{ultima}).

Incidentally, this shows that the open subsets $f_+(U_{\mu})$ yield a quasifold atlas for ${M}_{\Phi_Y\geq \varepsilon}$.
\qed

Consider now the open symplectic quasifolds
$$M_{\Phi_Y>\varepsilon}=\{m\in M\;|\;\Phi_Y(m)>\varepsilon\}$$
and
$$M_{\Phi_Y<\varepsilon}=\{m\in M\;|\;\Phi_Y(m)<\varepsilon\}.$$
\begin{remark}\label{apertone}\rm{
As in standard symplectic cutting, $M_{\Phi_Y>\varepsilon}$ identifies naturally with 
$$\{[m\!:\!w]\in \overline{M}_{\Phi_Y\geq \varepsilon}\;|\; |w|\neq0\}$$ and is therefore
open and dense in
$\overline{M}_{\Phi_Y\geq \varepsilon}$. Moreover, the difference 
$\overline{M}_{\Phi_Y\geq \varepsilon}\setminus M_{\Phi_Y>\varepsilon}$ can be identified with the orbit space
$\Phi_Y^{-1}(\varepsilon)/D^1$. Similarly, $M_{\Phi_Y<\varepsilon}$
is open and dense in
$\overline{M}_{\Phi_Y\leq \varepsilon}$ and the difference 
$\overline{M}_{\Phi_Y\leq \varepsilon}\setminus M_{\Phi_Y<\varepsilon}$ can be also identified with
$\Phi_Y^{-1}(\varepsilon)/D^1$.
}
\end{remark}
Notice that, under the assumptions of Proposition~\ref{sottoquasireticolo}, if  $\widetilde{\Phi}_Y$ denotes the component 
$\langle \widetilde{\Phi},Y\rangle$, we have the following
\begin{prop}\label{boh}
$$M_{\Phi_Y>\varepsilon}\simeq\left(\widetilde{M}_{\widetilde{\Phi}_Y>\varepsilon}\right)/\Gamma
\quad\hbox{and}\quad M_{\Phi_Y<\varepsilon}\simeq\left(\widetilde{M}_{\widetilde{\Phi}_Y<\varepsilon}\right)/\Gamma.
$$
\end{prop}
\begin{remark}\label{liscio}{\rm
It is easily shown that $\Phi_Y$ establishes a one--to--one correspondence between the fixed points of the $D^1$--action on 
$\Phi_Y^{-1}(\varepsilon)=\Phi^{-1}(\Delta\cap H(Y,\varepsilon))$ and the vertices of $\Delta$ contained in $H(Y,\varepsilon)$.
Recall that, in the classical cutting procedure, $M$ is a manifold and $D^1$ is a circle $S^1$ which is assumed to act freely
on $\Phi_Y^{-1}(\varepsilon)$; therefore the quotient $\Phi_Y^{-1}(\varepsilon)/S^1$ is a smooth
manifold, and, if $Y$ is primitive, the symplectic cuts are also smooth manifolds. Notice that,
if $S^1$ acts freely on $\Phi_Y^{-1}(\varepsilon)$, $H(Y,\varepsilon)$ does not contain vertices of $\Delta$.
However, in our setting we only require that the cut polyhedra are $n$--dimensional and simple.
In fact,  if this is the case, even if $H(Y,\varepsilon)$ contains some vertices of $\Delta$, the quotient
$\Phi_Y^{-1}(\varepsilon)/D^1$ and the symplectic cuts turn out to be honest quasifolds.
Incidentally, this applies also to the symplectic toric manifold case: even if the circle $S^1$ fixes some points in
the level set $\Phi_Y^{-1}(\varepsilon)$,  
when the cut polyhedra are Delzant, the quotient $\Phi_Y^{-1}(\varepsilon)/S^1$ and the symplectic cuts are smooth manifolds.
This provides
examples of Guillemin--Sternberg's {\em simple critical level sets} \cite{gs1}, meaning level sets 
of Hamiltonian circle actions that are singular but give rise to smooth reduced spaces. 
Also, under the standard assumptions, 
there is a reversing--orientation bundle--isomorphism between
the normal bundles of the submanifold 
$\Phi_Y^{-1}(\varepsilon)/S^1$ in $\overline{M}_{\Phi_Y\geq \varepsilon}$ and 
$\overline{M}_{\Phi_Y\leq \varepsilon}$ (see \cite[Section 1.1]{lerman}). 
This is not the case when the $S^1$--action fixes some points in $\Phi_Y^{-1}(\varepsilon)$. In fact, let us
cut the Delzant square $[0,1]\times[0,1]$ with the line $H(Y,1)$ containing
the vertices $(1,0)$ and $(0,1)$. The corresponding 
cut spaces are both equal to $\C P^2$ and, moreover, 
$\Phi_Y^{-1}(1)/S^1$ is equal to $\C P^1$; therefore the normal bundle of $\Phi_Y^{-1}(1)/S^1$ in both cut spaces is $O(1)$.
Similar remarks hold in the orbifold case.}
\end{remark}
\begin{remark}{\rm When $H(Y,\varepsilon)$ does not contain vertices of $\Delta$, the cut spaces and the quotient $\Phi_Y^{-1}(\varepsilon)/D^1$ are examples of reduced spaces for symplectic toric quasifolds, as defined in \cite{reduction}.}
\end{remark}
\section{Cutting the Penrose kite}
\label{cutkite}
\begin{figure}[h]
\begin{center}
\includegraphics[width=10cm]{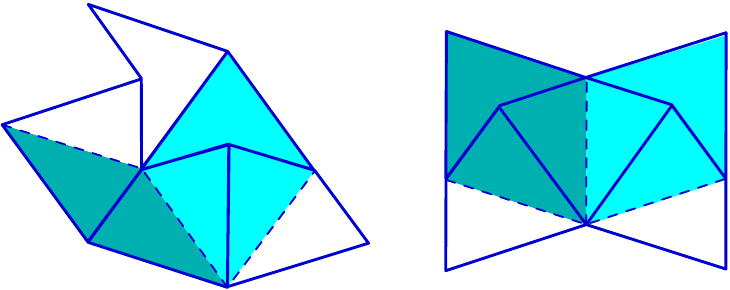}
\caption{From kite and dart to rhombus tilings}
\label{passaggio}
\end{center}
\end{figure}

Penrose tilings are a good source of significant examples of symplectic toric quasifolds \cite{rhombus,kite}. The Penrose kite, for example, 
is the simplest nonrational polytope; we show in \cite[Theorem~6.1]{kite} that the corresponding quasifold is not a global quotient of a manifold
modulo the action of a discrete group. As it turns out, rhombus tilings can be obtained from kite and dart tilings by cutting all 
kites in half (see Figure~\ref{passaggio}). 
It is exactly this picture that inspired us to consider symplectic
cuts in the nonrational setting, including the case where the cutting hyperplane meets the vertices of $\Delta$
(cf. Remark~\ref{liscio}).

\begin{figure}[h]
\begin{minipage}[b]{0.45\linewidth}
\begin{center}
\includegraphics[width=5cm]{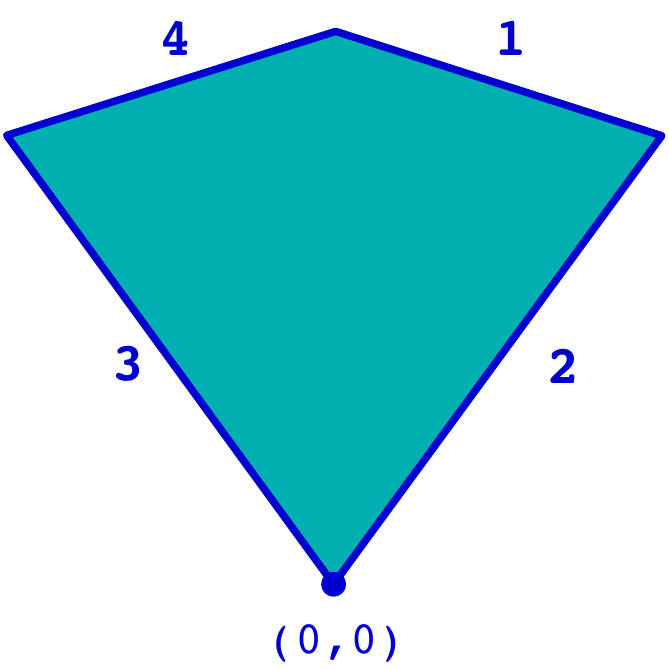}
\end{center}
\caption{The kite} \label{kite}
\end{minipage}
\qquad
\begin{minipage}[b]{0.45\linewidth}
\begin{center}
\includegraphics[width=4.5cm]{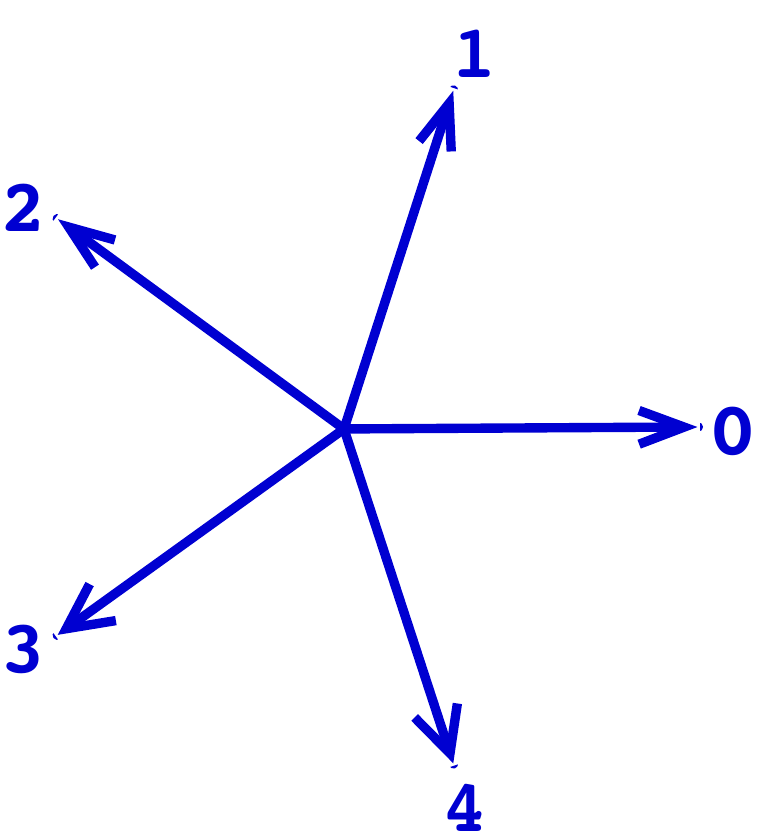}
\end{center}
\caption{The vectors $Y_0,Y_1,Y_2,Y_3,Y_4$}
\label{stargrid}
\end{minipage}
\end{figure}

\begin{figure}[h]
\begin{center}
\includegraphics[width=5cm]{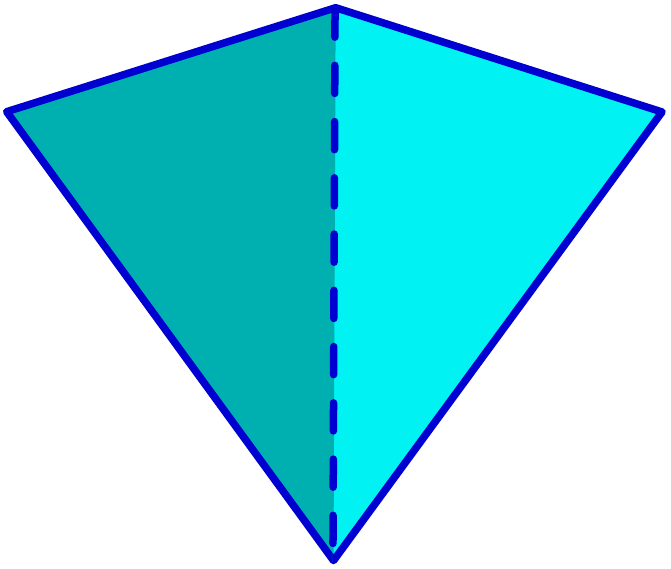}
\end{center}
\caption{The cut kite} \label{cutkitefig}
\end{figure}

Let us consider the Penrose kite $\Delta$ in Figure~\ref{kite} 
having 
vertices
$
(0,0)$,
$\left(\frac{1}{\phi},\frac{\phi}{\sqrt{2+\phi}}\right)$,
$\left(0,\frac{2}{\sqrt{2+\phi}}\right)$, and $\left(-\frac{1}{\phi},\frac{\phi}{\sqrt{2+\phi}}\right)$,
where $\phi=\frac{1+\sqrt{5}}{2}$ is the golden ratio; we recall that $\phi=\frac{1}{\phi}+1$.
The polytope $\Delta$ is quasirational with respect to the pentagonal quasilattice $Q$ in $\R^2$ spanned by the vectors
\begin{equation}\label{star}\left\{\begin{array}{l}
Y_0=(1,0)\\
Y_1=\frac{1}{2}(\frac{1}{\phi},\sqrt{2+\phi})\\
Y_2=\frac{1}{2}(-\phi,\frac{1}{\phi}\sqrt{2+\phi})\\
Y_3=\frac{1}{2}(-\phi,-\frac{1}{\phi}\sqrt{2+\phi})\\
Y_4=\frac{1}{2}(\frac{1}{\phi},-\sqrt{2+\phi})\\
\end{array}
\right.\end{equation}
(see
Figure~\ref{stargrid}).
In fact, $$\Delta=\cap_{j=1}^{4}\{\lambda\in(\R^2)^*\;|\;\langle \lambda,X_j\rangle\geq \lambda_j\},$$
where $X_1=-Y_1$, $X_2=Y_2$, $X_3=-Y_3$, $X_4=Y_4$,
$\lambda_1=\lambda_4=
-1$, and $\lambda_2=\lambda_3=0$.
By applying the generalized Delzant procedure (see Theorem~\ref{thmp1}) to $\D$ we obtain
\begin{equation}\label{quasikite}
M=\frac{\{(z_1,z_2,z_3,z_4)\in\C^4\;|\;\phi|z_1|^2+|z_2|^2+\phi|z_3|^2=\phi,\,-|z_1|^2+|z_2|^2+\phi|z_4|^2=\phi-1\}}
{\left\{\exp_{T^4}\left(s-(\phi-1)t,(\phi-1)(s+t),s,t\right)\;|\;
s,t \in\R\right\}}.
\end{equation}
We now cut the kite $\Delta$ with the line $H(Y_0,0)=\{\lambda\in(\R^2)^*\,|\,\langle\lambda,Y_0\rangle=0\}$
(see Figure~\ref{cutkitefig}). 
Notice that $\Lambda=\exp_{S^1}(\phi\Z)$. Consider the actions of $D^1$ on $M\times\C/\Lambda$:
$$
\exp_{D^1}(t)\cdot([z_1:z_2:z_3:z_4],[w])=\left(\left[e^{-2\pi i\frac{t}{\phi}}z_1:e^{-2\pi i t}z_2:z_3:
z_4\right],\left[e^{\mp 2\pi i t}w\right]\right).
$$
The corresponding moment mappings are:
$$\begin{array}{cll}\nu_{\mp} (([\vz],[w]))&=&\Phi_Y([\vz])\mp|w|^2\\
&=&-\frac{1}{\phi}(|z_1|^2-1)-|z_2|^2\mp|w|^2.
\end{array}
$$
Thus, by definition, the cut spaces are given by:
$$\overline{M}_{\Phi_Y\geq0}=\left\{([\vz],[w])\in M\times \C/\Lambda\;|\;
|z_1|^2+\phi|z_2|^2+\phi|w|^2=1\right\}/D^1
$$
and
$$\overline{M}_{\Phi_Y\leq0}=\left\{([\vz],[w])\in M\times \C/\Lambda\;|\;
\phi|z_3|^2+|z_4|^2+\phi|w|^2=1\right\}/D^1.$$
By Theorem~\ref{teorema},
they can be identified with the following quasifold, which corresponds, by the generalized Delzant procedure, to both the right and left triangle
$$M_{0^+}=M_{0^-}=\frac{\left\{ \,(u_1,u_2,u_3)\in\C^3\,|\,|u_1|^2+\phi|u_2|^2+\phi |u_3|^2=1
\,\right\}}{\left\{\,\exp\left((\phi-1)s,s+k\phi,s\right)\in T^3\,|\,s\in\R,k\in \Z\,\right\}}.$$

\section{Nonrational symplectic blow--ups}
\label{blowup}
Motivated by what happens in the classical setting \cite{lerman}, we define symplectic blowing up in the 
nonrational toric setting to be a special case of symplectic cutting in the following sense.
Let us consider an $n$--dimensional simple pointed polyhedron $\Delta\subset\rndu$ that is quasirational with respect to a quasilattice $Q$.
Then $\Delta=\cap_{j=1}^{d}\{\lambda\in\rndu\;|\;\langle \lambda,X_j\rangle\geq \lambda_j\}$, with
inward pointing normal vectors $\{X_1,\ldots,X_d\}$ in $Q$.
Let us apply the generalized Delzant construction to $\Delta$ with respect to $\{X_1,\ldots,X_d\}$ 
and let $M$ be the corresponding symplectic quasifold. Let $p\in M$ be a fixed point for the $D^n$--action and let $\nu=\Phi(p)$ be the corresponding vertex.
We can assume, up to a translation of $\Delta$, that $\nu=0$.
A blow--up of $M$ of an $\varepsilon$--amount at the fixed point $p$ is the symplectic cut $\overline{M}_{\Phi_Y\geq \varepsilon}$, where $Y\in Q$ and
$\epsilon\in\R$ are chosen so that $\D_{\varepsilon^-}$ has the combinatorial type of a simplex.
By Remark~\ref{apertone}, this is consistent with McDuff's approach to symplectic blow--ups \cite{md}, as already noticed, in the standard setting, in [10].
\begin{figure}[h]
\begin{center}
\includegraphics[width=6cm]{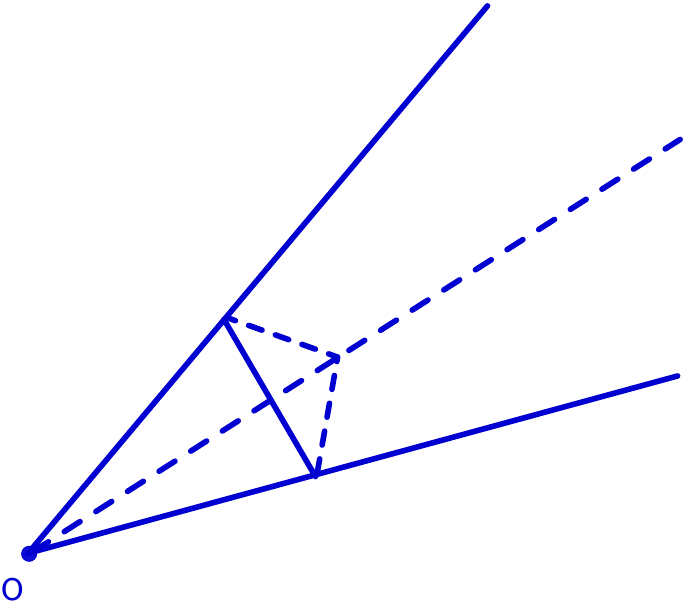}
\end{center}
\caption{The simple convex polyhedral cone ${\cal C}$} \label{cono}
\end{figure}
The following example provides the local model for blow--ups:
\begin{example}\label{local blow up}{\rm
We consider a simple $n$--dimensional convex polyhedral cone $\cal C$ in $(\R^n)^*$ with apex at the origin, and we assume that it is
quasirational with respect to a quasilattice $Q$ (see Figure~\ref{cono}). Then 
$${\cal C}=\cap_{j=1}^{n}\{\lambda\in\rndu\;|\;\langle \lambda,X_j\rangle\geq 0\},$$ with inward pointing normal vectors
$\{X_1,\ldots,X_n\}$ in $Q$. Let us apply the generalized Delzant procedure. In this case the projection 
$$\begin{array}{ccccc}
\pi&\colon&\R^n&\longrightarrow&\R^n\\
&&e_j&\longmapsto&X_j
\end{array}
$$
induces a group epimorphism 
$$\Pi \,\colon T^n=\R^n/\Z^n\longrightarrow D^n=\R^n/Q.$$ Here $N=\hbox{Ker}(\Pi)$ is the discrete group
$$\Gamma=\left\{(e^{2\pi i x_1},\ldots,e^{2\pi i x_n})\in T^n\;\Bigg|\;\sum_{j=1}^n x_j X_j\in Q\right\}.$$ Recall that
$D^n\simeq T^n/\Gamma$. 
The moment mapping with respect to the induced $D^n$--action on $\C^n/\Gamma$ is given by
$$\Phi([\vz])=\sum_{j=1}^n\alpha_j|z_j|^2,$$ where 
$\alpha_j=(\pi^*)^{-1}(e^*_j)$, $j=1,\ldots,n$, is the basis of $(\R^n)^*$ dual to the basis $\{X_1,\ldots,X_n\}$. 
The image of $\Phi$ is the cone $\cal C$ itself.
Let $$\hbox{int}({\cal C}^*)=\{X\in\R^n\;|\;\langle\alpha_j,X\rangle>0, j=1,\dots,n\}$$ and let 
$Y\in\hbox{int}({\cal C}^*)\cap Q$. 
Let us now cut the cone $\cal C$ with the hyperplane  
$H(Y,\varepsilon)$.
Since $Y=\sum_{j=1}^n\langle\alpha_j,Y\rangle X_j$, we have that
$$\Lambda=\{e^{2\pi i \tau}\in S^1\,|\,\tau Y\in Q\}=\left\{e^{2\pi i \tau}\in S^1\;\Bigg|\;\sum_{j=1}^n\tau\langle\alpha_j,Y\rangle X_j\in Q\right\}$$
and that the actions of $D^1$ on 
$(\C^n/\Gamma)\times(\C/\Lambda)$ are given by 
$$
\exp_{D^1}(t)\cdot([z_1:\cdots:z_n],[w])=\left(\left[e^{2\pi it\langle\alpha_1,Y\rangle}z_1:
\cdots:e^{2\pi it\langle\alpha_n,Y\rangle}z_n\right],\left[e^{\mp 2\pi i t}w\right]\right).
$$
Then the cut spaces are
$$(\overline{\C^n/\Gamma})_{\Phi_Y\geq\varepsilon}=\left\{([\vz],[w])\in (\C^n/\Gamma)\times 
(\C/\Lambda)\;\Bigg|\;\sum_{j=1}^n\langle\alpha_j,Y\rangle|z_j|^2-|w|^2=\varepsilon\right\}/D^1$$
and
$$(\overline{\C^n/\Gamma})_{\Phi_Y\leq\varepsilon}=\left\{([\vz],[w])\in (\C^n/\Gamma)\times 
(\C/\Lambda)\;\Bigg|\;\sum_{j=1}^n\langle\alpha_j,Y\rangle|z_j|^2+|w|^2=\varepsilon\right\}/D^1.$$
In particular, the quasifold $(\overline{\C^n/\Gamma})_{\Phi_Y\geq\varepsilon}$ is the symplectic blow--up of 
$\C^n/\Gamma$ at the origin of an $\varepsilon$--amount. 
By Theorem~\ref{teorema}, $(\overline{\C^n/\Gamma})_{\Phi_Y\geq\varepsilon}$ and $(\overline{\C^n/\Gamma})_{\Phi_Y\leq\varepsilon}$
can be identified with the symplectic toric quasifolds
$$(\C^n/\Gamma)_{\varepsilon^\pm}=\left\{\vu\in\C^{n+1}\;\Bigg|\;\sum_{j=1}^n\langle\alpha_j,Y\rangle|u_j|^2\mp|u_{n+1}|^2
=\varepsilon\right\}/N_{\varepsilon^\pm}$$
where
$$N_{\varepsilon^\pm}=\{(\gamma_1 e^{2\pi i \theta\langle\alpha_1,Y\rangle},\ldots,\gamma_n e^{2\pi i \theta\langle\alpha_n,Y\rangle},e^{\mp 2\pi i \theta})
\in T^{n+1}\,|\,\gamma=(\gamma_1,\ldots,\gamma_n)\in\Gamma,\theta\in\R\}.$$ 
The key point here is to notice that $e^{2\pi\tau}\in\Lambda$ if, and only if, 
$(e^{2\pi i\tau\langle\alpha_1,Y\rangle},\ldots,e^{2\pi i\tau\langle\alpha_n,Y\rangle})\in\Gamma$.}
\end{example}
\begin{example}\label{quadrante}
{\rm As a special case of the previous example, let ${\cal C}$ be the positive quadrant in $(\R^2)^*$. Notice that
$${\cal C}=\{\lambda\in(\R^2)^*\;|\;\langle \lambda,e_1\rangle\geq 0\}\cap\{\lambda\in(\R^2)^*\;|\;\langle \lambda,e_2\rangle\geq 0\}.$$
Then, if $Q$ is a quasilattice containing $\Z^2$, the cone ${\cal C}$ is quasirational with respect to $Q$. 
In this case the quasitorus $D^2=\R^2/Q$ acts on $\C^2/\Gamma$ with
moment mapping $\Phi([\vz])=(|z_1|^2,|z_2|^2)$.  
\begin{figure}[h]
\begin{center}
\includegraphics[width=6cm]{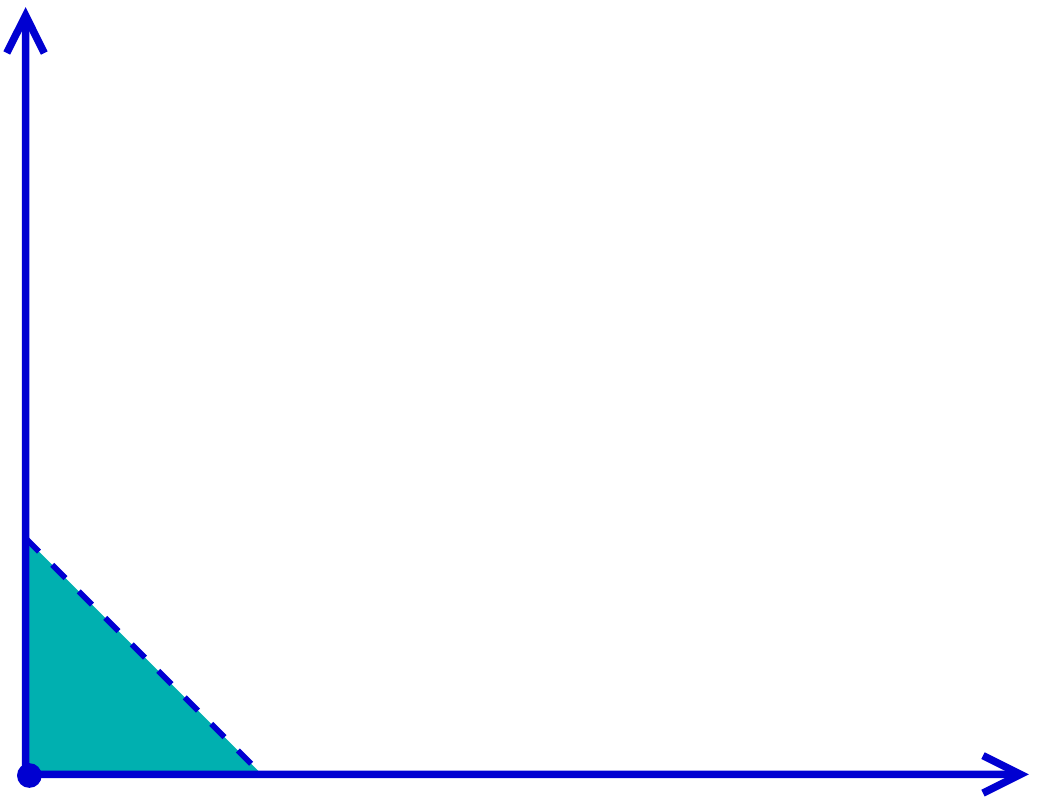}
\end{center}
\caption{Blowing up the origin in $\C^2/\Gamma$} \label{cutquadrant}
\end{figure}
Let us now cut the cone $\cal C$ with the line $H(Y,\varepsilon)$, with $Y=(1,1)$ and $\varepsilon>0$ (see Figure~\ref{cutquadrant}). 
We obtain the two cut spaces
$$(\overline{\C^2/\Gamma})_{\Phi_Y\geq\varepsilon}=\left\{([\vz],[w])\in(\C^2/\Gamma)\times(\C/\Lambda)\;\Big|\;
|z_1|^2+|z_2|^2=|w|^2+\varepsilon\right\}/D^1$$
and
$$(\overline{\C^2/\Gamma})_{\Phi_Y\leq\varepsilon}=\left\{([\vz],[w])\in
(\C^2/\Gamma)\times(\C/\Lambda)\;\Big|\;
|z_1|^2+|z_2|^2+|w|^2=\varepsilon\right\}/D^1.$$
The first quotient is the blow--up at the origin of $\C^2/\Gamma$ of an $\varepsilon$--amount. 
The quasisphere $$S^2/\Lambda=\{|z_1|^2+|z_2|^2=\varepsilon\}/D^1$$ is naturally embedded in both spaces as the set of
classes with $w=0$.  
By Theorem~\ref{teorema}, $(\overline{\C^2/\Gamma})_{\Phi_Y\geq\varepsilon}$ and 
$(\overline{\C^2/\Gamma})_{\Phi_Y\leq\varepsilon}$ can be identified 
with the symplectic toric quasifolds
$$(\C^2/\Gamma)_{\varepsilon^\pm}=\left\{\vu\in\C^{3}\,|\,|u_1|^2+|u_2|^2\mp|u_{3}|^2=\varepsilon\right\}/N_{\varepsilon^\pm}$$
where
$$N_{\varepsilon^\pm}=\{(\gamma_1 e^{2\pi i \theta},\gamma_2 e^{2\pi i \theta},e^{\mp 2\pi i \theta})\in T^3\,|\,
\gamma=(\gamma_1,\gamma_2)\in\Gamma,\theta\in\R\}.$$
From this one deduces, in particular, that $(\overline{\C^2/\Gamma})_{\Phi_Y\leq 1}\simeq\C P^2/\Gamma$.}
\end{example}

\section{Cutting and blowing up in arbitrary directions}
\label{arbitrary}
Let $\widetilde Q$  be a quasilattice in $\R^n$ and let $\Delta\in (\R^n)^*$ be an $n$--dimensional simple pointed polyhedron
that is quasirational with respect to $\widetilde Q$. 
Consider the symplectic toric quasifold $\widetilde M$ corresponding to 
$(\Delta,{\widetilde Q},\{X_1,\ldots,X_d\})$. We would like to be able to cut $\Delta$ with a hyperplane 
$H(Y,\varepsilon)$,
where $Y\in\R^n$ is no longer necessarily in $\widetilde Q$; 
the only remaining requirement is that the cut polytopes are $n$--dimensional and simple.
The idea is to consider the new quasilattice $Q=\hbox{Span}_{\Z}({\widetilde Q}\cup Y)$ and to apply the
generalized Delzant procedure to $(\Delta,Q,\{X_1,\ldots,X_d\})$. Then we
cut the corresponding quasifold $M$; this is now allowed since $Y\in Q$. 
By Proposition~\ref{sottoquasireticolo}, we have that $M\simeq{\widetilde M}/\Gamma$. Moreover, by Proposition~\ref{boh},
$\left(\widetilde{M}_{\widetilde{\Phi}_Y>\varepsilon}\right)/\Gamma$
embeds as an open dense
subset of $\overline{M}_{\Phi_Y\geq \varepsilon}$, while
$\left(\widetilde{M}_{\widetilde{\Phi}_Y<\varepsilon}\right)/\Gamma$
embeds as an open dense
subset of $\overline{M}_{\Phi_Y\leq \varepsilon}$.
It is therefore natural to define
$\overline{M}_{\Phi_Y\geq\varepsilon}$ and $\overline{M}_{\Phi_Y\leq\varepsilon}$ 
as cut spaces for the quasifold ${\widetilde M}$ as well.
We conclude by applying this procedure to a simple, yet representative example. 
For another application of this procedure, we refer the reader to \cite{hirze}.
\begin{example}[Cutting $\C^2$ in an arbitrary direction]{\rm Let us consider again the positive quadrant $\cal C$
of Example~\ref{quadrante}. The cone $\cal C$ is obviously rational with respect to the lattice $\Z^2$ and the symplectic
toric manifold corresponding to $({\cal C},\Z^2,\{e_1,e_2\})$ is $\C^2$ with the standard Hamiltonian action of $T^2=\R^2/\Z^2$.
\begin{figure}[h]
\begin{center}
\includegraphics[width=6cm]{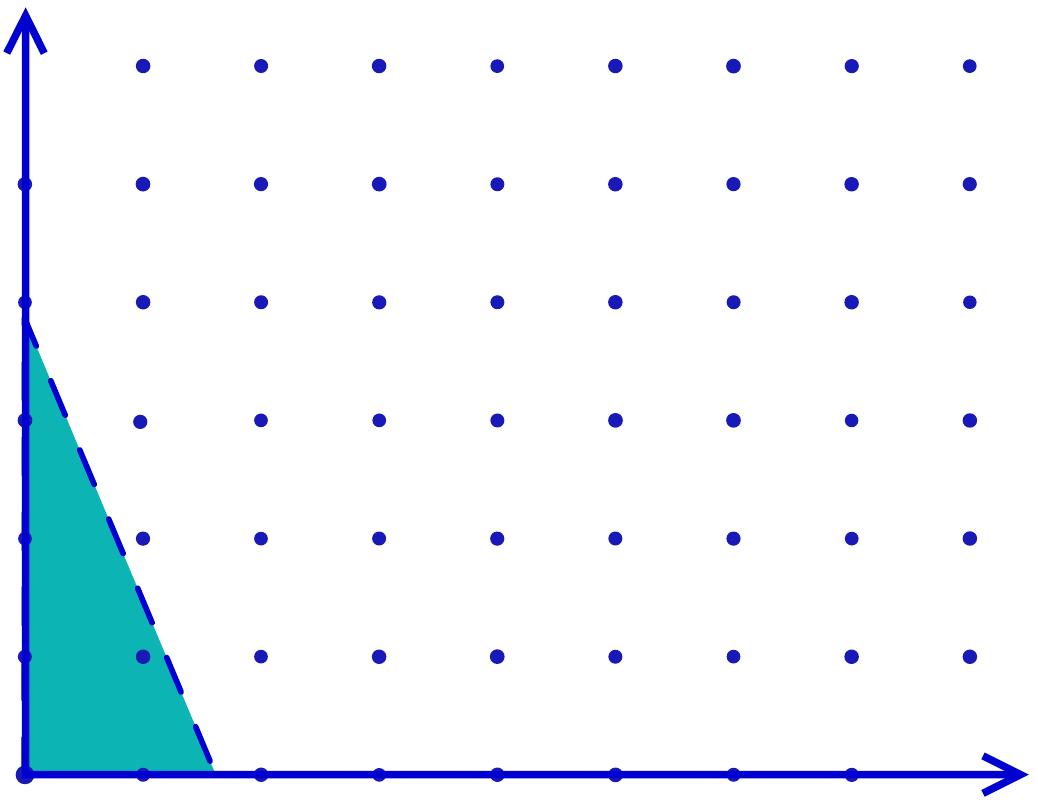}
\end{center}
\caption{An arbitrary cut of $\C^2$} \label{cutquadrant2}
\end{figure}

Consider now two positive real numbers $s,t$ such that  $\frac{s}{t}\notin \Q$. We cut the positive quadrant 
$\cal C$ with the line $H(Y,\varepsilon)$, where $Y=(s,t)$ and $\varepsilon>0$
(see Figure~\ref{cutquadrant2}). 
Let us consider the quasilattice $Q=\hbox{Span}_{\Z}(\Z^2\cup Y)$ and let us
apply the generalized Delzant procedure to $({\cal C},Q,\{e_1,e_2\})$. The corresponding quasifold is $M=\C^2/\Gamma$, 
where the action of $\Gamma=\exp_{T^2}(Q)$ on $\C^2$ is given by
$$(e^{2\pi i ls},e^{2\pi i lt})\cdot(z_1,z_2)=(e^{2\pi i ls}z_1,e^{2\pi i lt}z_2),\quad l\in\Z.$$
Since $\frac{s}{t}\notin\Q$ it is easily verified that $\Lambda$ is trivial
and that $D^1\simeq S^1$ acts on 
$(\C^2/\Gamma)\times\C$ as follows 
$$
\exp_{D^1}(\sigma)\cdot([z_1:z_2],w)=\left(\left[e^{2\pi i s\sigma}z_1:
e^{2\pi i t\sigma}z_2\right],e^{\mp 2\pi i \sigma}w\right).
$$
We obtain therefore the two cut spaces:
$$(\overline{\C^2/\Gamma})_{\Phi_Y\geq\varepsilon}=\left\{([\vz],w)\in(\C^2/\Gamma)\times \C\;\Big|\;
s|z_1|^2+t|z_2|^2=|w|^2+\varepsilon\right\}/D^1$$
$$(\overline{\C^2/\Gamma})_{\Phi_Y\leq\varepsilon}=\left\{([\vz],w)\in(\C^2/\Gamma)\times \C\;\Big|\;
s|z_1|^2+t|z_2|^2+|w|^2=\varepsilon\right\}/D^1.$$
The open subsets
$$\{(z_1,z_2)\in\C^2\;|\;s|z_1|^2+t|z_2|^2>\varepsilon\}/\Gamma$$ and $$\{(z_1,z_2)\in\C^2\;|\;s|z_1|^2+t|z_2|^2<\varepsilon\}/\Gamma$$ 
embed as open subsets in $(\overline{\C^2/\Gamma})_{\Phi_Y\geq\varepsilon}$ and $(\overline{\C^2/\Gamma})_{\Phi_Y\leq\varepsilon}$ respectively. 
By Theorem~\ref{teorema} $(\overline{\C^2/\Gamma})_{\Phi_Y\geq\varepsilon}$ and $(\overline{\C^2/\Gamma})_{\Phi_Y\leq\varepsilon}$ 
can be identified with the symplectic toric quasifolds
$$(\C^2/\Gamma)_{\varepsilon^\pm}=\frac{\left\{\vu\in\C^{3}\,|\,s|u_1|^2+t|u_2|^2\mp|u_{3}|^2=\varepsilon\right\}}
{\{(e^{2\pi i s\theta},e^{2\pi i t\theta},e^{\mp 2\pi i \theta})\in T^3\,|\,\theta\in\R\}}.$$
The quotient $(\C^2/\Gamma)_{\varepsilon^-}$ is
the natural generalization, in the quasifold setting, of the complex projective space $\C P^2$ (it is the {\em projective quasispace} of \cite[Example~3.6]{p}).}
\end{example}
\section*{Acknowledgements}
We thank Dan Zaffran for drawing our attention to the interest of blowing up in arbitrary directions \cite{bz}.
This research was partially supported by  grant  PRIN 2010NNBZ78\_$\!$\_012  (MIUR, Italy).

\noindent \sc Dipartimento di Matematica e Informatica "U. Dini",
Universit\`a di Firenze, Via S. Marta 3, 50139 Firenze,
ITALY, {\tt fiammetta.battaglia@unifi.it}\\
                        and\\
                        Dipartimento di Matematica e Informatica "U. Dini", Universit\`a di Firenze,
                        Piazza Ghiberti 27, 50122 Firenze, ITALY, {\tt elisa.prato@unifi.it}

\begin{thebibliography}{00}

\bibitem{a} M. Atiyah, Convexity and Commuting Hamiltonians,
{\em Bull. London Math. Soc.} \textbf{14} (1982), 1--15.

\bibitem{cx} F. Battaglia, E. Prato, Generalized toric varieties for
simple nonrational convex polytopes, Intern. Math. Res. Notices 24 (2001), 1315--1337.

\bibitem{rhombus} F. Battaglia, E. Prato, The Symplectic Geometry of Penrose Rhombus Tilings, 
\textit{J. Symplectic Geom.} \textbf{6} (2008), 139--158.

\bibitem{kite} F. Battaglia, E. Prato, The symplectic Penrose kite, {\em Comm. Math. Phys.} \textbf{299} (2010), 577--601.

\bibitem{reduction} F. Battaglia, E. Prato, Nonrational symplectic toric reduction, arXiv:1806.10431 [math.SG] (2018).

\bibitem{hirze} F. Battaglia, E. Prato, D. Zaffran, Hirzebruch surfaces in a one--parameter family, arXiv:1804.08503 [math.SG] (2018).

\bibitem{bz} F. Battaglia, D. Zaffran, LVMB manifolds as equivariant group compatifications, {\em in preparation}.

\bibitem{d} T. Delzant, Hamiltoniens p\'eriodiques et images convexes
de l'application moment, {\em Bull. S.M.F.} \textbf{116} (1988),
315--339.

\bibitem{gromov} M. Gromov, Partial Differential Relations, {\em Ergebnisse der Math.} \textbf{9}, Springer-Verlag (1986).

\bibitem{gs} V. Guillemin and S. Sternberg, Convexity Properties
of the Moment Mapping, {\em Invent. Math.} \textbf{67} (1982),
491--513.

\bibitem{gs1} V. Guillemin, S. Sternberg, Birational equivalence in the symplectic
category, {\em Invent. Math.}  \textbf{97} (1989),
485--522.

\bibitem{lerman} E. Lerman, Symplectic cuts, {\em Math. Res. Lett.} \textbf{2} (1995), 247--258.

\bibitem{lertol} E. Lerman, S. Tolman, Hamiltonian torus actions on
symplectic orbifolds and toric varieties, {\em Trans. Amer. Math. Soc.} \textbf{349} (1997), 4201--4230.

\bibitem{md0} D. McDuff, Examples of simply--connected symplectic non--K\"ahlerian manifolds, J. Differ. Geom.
\textbf{20} (1984) 267--277.

\bibitem{md} D. McDuff, Blow--ups and symplectic embeddings in dimension $4$, {\em Topology} \textbf{30} (1991), 409--421.

\bibitem{pcras} E. Prato, Sur une g\'en\'eralisation de la notion de V--vari\'et\'e, {\em C. R. Acad. Sci. Paris, 
Ser. I} \textbf{328} (1999), 887--890.

\bibitem{p}  E. Prato, Simple non--rational convex
polytopes via symplectic geometry, {\em Topology} \textbf{40}
(2001), 961--975.

\bibitem{ziegler} G. Ziegler, Lectures on polytopes, {\em Graduate Texts in Mathematics} \textbf{52}, Springer--Verlag (1995).

\end{thebibliography}
\end{document}